\documentclass[12pt,reqno]{amsart}
\usepackage{times}

\newtheorem{thm}{Theorem}
\newtheorem{lemma}[thm]{Lemma}
\newtheorem{cor}[thm]{Corollary}

\newtheorem*{mthm}{Main Theorem}

\newcommand{\C}{{\mathbb C}}
\newcommand{\cn}{{\mathbb C}^n}
\newcommand{\inc}{\int_{\C}}
\newcommand{\incn}{\int_{\cn}}
\newcommand{\dva}{\,dv_t}

\begin{document}

\title[Integral Operators]
{Integral Operators Induced by the Fock Kernel}

\author{ Milutin Dostani\'c and Kehe Zhu}
\address{Matemati\v cki fakultet\\ Studentski trg 16\\ 11000 Beograd\\ Serbia}
\email{domi@matf.bg.ac.yu}
\address{Department of Mathematics\\ SUNY\\ Albany, NY 12222, USA}
\email{kzhu@math.albany.edu }
\date{December 20, 2006}
\thanks{Dostani\'c is partialy supported by MNZZS Grant $N^o ON 144010$}
\thanks{Zhu is partially supported by the US National Science Foundation}

\begin{abstract}
We study the $L^p$ boundedness  and find the norm of a class of
integral operators induced by the reproducing kernel of Fock
spaces over $\cn$.
\end{abstract}

\maketitle

\section{Introduction}

Our analysis will take place in the $n$-dimensional complex Euclidean space
$\cn$. For any two points $z=(z_1,\cdots,z_n)$ and $w=(w_1,\cdots,w_n)$ in
$\cn$ we write
$$\langle z,w\rangle=z_1\overline w_1+\cdots+z_n\overline w_n,$$
and
$$|z|=\sqrt{|z_1|^2+\cdots+|z_n|^2}.$$

For any $t>0$ we consider the Gaussian probability measure
$$dv_t(z)=\left(\frac t\pi\right)^ne^{-t|z|^2}\,dv(z)$$
on $\cn$, where $dv$ is ordinary Lebesgue volume measure on $\cn$. Let
$H(\cn)$ denote the space of all entire functions on $\cn$. We then define
$$F^p_t=L^p(\cn,dv_t)\cap H(\cn)$$
for $0<p<\infty$. These spaces are often called Fock spaces, or Segal-Bargman
spaces, over $\cn$. See \cite{bargmann}\cite{BC1}\cite{BC2}
\cite{folland}\cite{guillemin}\cite{JPR}\cite{T}\cite{XZ}.

For $p>0$ and $t>0$ we are going to write
$$\|f\|_{t,p}=\left[\incn|f(z)|^p\dva(z)\right]^{\frac1p},$$
and
$$\langle f,g\rangle_t=\incn f(z)\,\overline{g(z)}\,dv_t(z).$$

It is well known that each Fock space $F^p_t$ is a closed linear subspace
of $L^p(\cn,dv_t)$. In particular, in the Hilbert space setting of
$L^2$, there exists a unique orthogonal projection
$$P_t:L^2(\cn,dv_t)\to F^2_t.$$
Furthermore, this projection coincides with the restriction of the
following integral operator to $L^2(\cn,dv_t)$:
\begin{equation}
S_t f(z)=\incn e^{t\langle z,w\rangle}f(w)\,dv_t(w).
\label{eq1}
\end{equation}
The integral kernel above,
$$K_t(z,w)=e^{t\langle z,w\rangle},$$
is the reproducing kernel of $F^2_t$.

The purpose of this paper is to study the action of the operator $S_t$
on the spaces $L^p(\cn,dv_s)$, where $s>0$. We also consider the closely
related integral operator
$$T_t f(z)=\incn|K_t(z,w)|f(w)\,dv_t(w),$$
or more explicitly,
\begin{equation}
T_tf(z)=\incn|e^{t\langle z,w\rangle}|f(w)\,dv_t(w).
\label{eq2}
\end{equation}

The main result of the paper is the following.

\begin{mthm}
Suppose $t>0$, $s>0$, and $p\ge1$. Then
the following conditions are equivalent.
\begin{enumerate}
\item[(a)] $T_t$ is bounded on $L^p(\cn,dv_s)$.
\item[(b)] $S_t$ is bounded on $L^p(\cn,dv_s)$.
\item[(c)] $pt=2s$.
\end{enumerate}
Furthermore, the norms of $T_t$ and $S_t$ on $L^p(\cn,dv_s)$ satisfy
$$\|S_t\|\le\|T_t\|= 2^n$$
whenever $pt=2s$.
\end{mthm}

The equivalence of (a), (b), and (c) is not new; it is implicit in \cite{JPR}
for example. So our main contribution here is the identity $\|T_t\|=2^n$. The
accurate calculation of the norm of an integral operator is an interesting but often 
difficult problem. We mention a few successful examples in the literature: the norm 
of the Cauchy projection on $L^p$ of the unit circle is determined in \cite{pic}, the 
norm of the Cauchy projection on $L^p$ spaces of more general domains is estimated
in \cite{dos1}, an asymptotic formula for the norm of the Bergman projection on 
$L^p$ spaces of the unit ball is given in \cite{zhu3}, and the norm of the Berezin 
transform on the unit disk is calculated in \cite{dos2}.

As a consequence of the theorem above, we see that the densely defined operator
$$P_t:L^p(\cn,dv_t)\to L^p(\cn,dv_t)$$
is unbounded for any $p\not=2$. This is in sharp contrast to the theory of Hardy 
spaces and the theory of Bergman spaces. For example, if $P$ is the Bergman projection 
for the open unit ball $B_n$, that is, if $P$ is the orthogonal projection
$$P: L^2(B_n,dv)\to L^2(B_n,dv)\cap H(B_n),$$
where $H(B_n)$ is the space of holomorphic functions in $B_n$, then
$$P: L^p(B_n,dv)\to L^p(B_n,dv)$$
is bounded for every $p>1$. A similar result holds for the Cauchy-Sz\"ego
projection in the theory of Hardy spaces. See \cite{rudin} and \cite{zhu2}.

A more general class of integral operators induced by the Bergman kernel
on the unit ball $B_n$ have been studied in \cite{FR}\cite{KZ}\cite{zhu1}.

We wish to thank Peter Duren and James Tung for bringing to our attention the
references \cite{JPR} and \cite{S}.

\section{Preliminaries}

For an $n$-tuple $m=(m_1,\cdots,m_n)$ of nonnegative integers we
are going to write
$$|m|=m_1+\cdots+m_n,\quad m!=m_1!\cdots m_n!.$$
If $z\in\cn$, we also write
$$z^m=z_1^{m_1}\cdots z_n^{m_n}.$$
When the dimension $n$ is 1, we use $dA$ instead of $dv$, and $dA_t$
instead of $dv_t$. Thus for $t>0$ and $z\in\C$, we have
$$dA_t(z)=\frac t\pi \,e^{-t|z|^2}\,dA(z),$$
where $dA$ is ordinary area measure on the complex plane $\C$.

\begin{lemma}
Let $m=(m_1,\cdots,m_n)$ be an $n$-tuple of nonnegative integers.
For any $t>0$ and $p>0$ we have
$$\incn|z^m|^p\,dv_t(z)=\prod_{k=1}^n\frac{\Gamma((pm_k/2)+1)}
{t^{pm_k/2}}.$$
In particular,
$$\incn|z^m|^2\dva(z)=\frac{m!}{t^{|m|}}.$$
\label{1}
\end{lemma}

\begin{proof}
We evaluate the integral in polar coordinates.
\begin{eqnarray*}
\incn|z^m|^p\,dv_t(z)&=&\prod_{k=1}^n\inc|z_k|^{pm_k}\,dA_t(z_k)\\
&=&\prod_{k=1}^n\frac t\pi\inc|z_k|^{pm_k}e^{-t|z_k|^2}\,dA(z_k)\\
&=&\prod_{k=1}^n2t\int_0^\infty r^{pm_k+1}e^{-tr^2}\,dr\\
&=&\prod_{k=1}^nt\int_0^\infty r^{pm_k/2}e^{-tr}\,dr\\
&=&\prod_{k=1}^n\frac1{t^{pm_k/2}}\int_0^\infty r^{pm_k/2}e^{-r}\,dr\\
&=&\prod_{k=1}^n\frac{\Gamma((pm_k/2)+1)}{t^{pm_k/2}}.
\end{eqnarray*}
The second integral is obviously a special case of the first one.
\end{proof}

Recall that the restriction of the operator $S_t$ to $L^2(\cn,dv_t)$ is
nothing but the orthogonal projection onto $F^2_t$. Consequently, we have
the following reproducing formula.

\begin{lemma}
If $f$ is in $F^2_t$, then $S_tf=f$, that is,
$$f(a)=\incn e^{t\langle a,z\rangle}\,f(z)\,dv_t(z)$$
for all $a\in\cn$.
\label{2}
\end{lemma}

A special case of the reproducing formula above is the following:
\begin{equation}
K_t(a,a)=\incn|K_t(a,z)|^2\,dv_t(z),\qquad a\in\cn.
\label{eq3}
\end{equation}
As an application of this identity, we obtain the following fundamental
integrals for powers of kernel functions in Fock spaces.

\begin{lemma}
Suppose $t>0$ and $s$ is real. Then
$$\incn|e^{s\langle z,a\rangle}|\,dv_t(z)=e^{s^2|a|^2/4t}$$
for all $a\in\cn$.
\label{3}
\end{lemma}

\begin{proof}
It follows from (\ref{eq3}) that
\begin{eqnarray*}
\incn|e^{s\langle z,a\rangle}|\,dv_t(z)&=&\incn|e^{t\langle sa/2t,z
\rangle}|^2\,dv_t(z)\\
&=&\incn|K_t(sa/2t,z)|^2\,dv_t(z)\\
&=&K_t(sa/2t,sa/2t)\cr\noalign{\vskip 5pt}
&=&e^{s^2|a|^2/4t}.
\end{eqnarray*}
This proves the desired identity.
\end{proof}

We need two well-known results from the theory of integral operators.
The first one concerns the adjoint of a bounded integral operator.

\begin{lemma}
Suppose $1\le p<\infty$ and $1/p+1/q=1$. If an integral operator
$$Tf(x)=\int_X K(x,y)f(y)\,d\mu(y)$$
is bounded on $L^p(X,d\mu)$, then its adjoint
$$T^*:L^q(X,d\mu)\to L^q(X,d\mu)$$
is the integral operator given by
$$T^*f(x)=\int_X\overline{K(y,x)}f(y)\,d\mu(y).$$
\label{4}
\end{lemma}

\begin{proof}
See \cite{HS} for example.
\end{proof}

The second result is a useful criterion for the boundedness of integral
operators on $L^p$ spaces, usually referred to as Schur's test.

\begin{lemma}
Suppose $H(x,y)$ is a positive kernel and
$$Tf(x)=\int_X H(x,y)f(y)\,d\mu(y)$$
is the associated integral operator. Let $1<p<\infty$ with $1/p+1/q=1$. If
there exists a positive function $h(x)$ and positive constants $C_1$ and
$C_2$ such that
$$\int_XH(x,y)h(y)^q\,d\mu(y)\le C_1h(x)^q,\qquad x\in X,$$
and
$$\int_XH(x,y)h(x)^p\,d\mu(x)\le C_2h(y)^p,\qquad y\in X,$$
then the operator $T$ is bounded on $L^p(X,d\mu)$. Moreover, the norm
of $T$ on $L^p(X,d\mu)$ does not exceed $C_1^{1/q}C_2^{1/p}$.
\label{5}
\end{lemma}

\begin{proof}
See \cite{zhu2} for example.
\end{proof}

\section{Integral Operators Induced by the Fock Kernel}

For any $s>0$ we rewrite the integral operators $S_t$ and $T_t$ defined in
(\ref{eq1}) and (\ref{eq2}) as follows.
$$S_t f(z)=\left(\frac ts\right)^n\incn e^{t\langle z,w\rangle
+s|w|^2-t|w|^2}f(w)\,dv_s(w),$$
and
$$T_t f(z)=\left(\frac ts\right)^n\incn|e^{t\langle z,w\rangle
+s|w|^2-t|w|^2}|\,f(w)\,dv_s(w).$$
It follows from Lemma~\ref{4} that the adjoint of $S_t$ and $T_t$ with respect
to the integral pairing
$$\langle f,g\rangle_s=\incn f(z)\overline{g(z)}\,dv_s(z)$$
is given respectively by
\begin{equation}
S_t^*f(z)=\left(\frac ts\right)^ne^{(s-t)|z|^2}\incn
e^{t\langle z,w\rangle}f(w)\,dv_s(w),
\label{eq4}
\end{equation}
and
\begin{equation}
T_t^*f(z)=\left(\frac ts\right)^ne^{(s-t)|z|^2}\incn
|e^{t\langle z,w\rangle}|f(w)\,dv_s(w).
\label{eq5}
\end{equation}

We first prove several necessary conditions for the operator $S_t$
to be bounded on $L^p(\cn,dv_s)$.

\begin{lemma}
Suppose $0<p<\infty$, $t>0$, and $s>0$. If $S_t$ is bounded
on $L^p(\cn,dv_s)$, then $pt\le 2s$.
\label{6}
\end{lemma}

\begin{proof}
Consider functions of the following form:
$$f_{x,k}(z)=e^{-x|z|^2}z_1^k,\qquad z\in\cn,$$
where $x>0$ and $k$ is a positive integer.

We first use Lemma~\ref{1} to calculate the norm of $f_{x,k}$ in
$L^p(\cn,dv_s)$.
\begin{eqnarray*}
\incn|f_{x,k}|^p\,dv_s&=&\left(\frac s\pi\right)^n\incn|z_1|^{pk}
e^{-(px+s)|z|^2}\,dv(z)\\
&=&\left(\frac s{px+s}\right)^n\incn|z_1^k|^p\,dv_{px+s}(z)\\
&=&\left(\frac s{px+s}\right)^n\frac{\Gamma((pk/2)+1)}{(px+s)^{pk/2}}.
\end{eqnarray*}

We then calculate the closed form of $S_t(f_{x,k})$ using the
reproducing formula from Lemma~\ref{2}.
\begin{eqnarray*}
S_t(f_{x,k})(z)&=&\left(\frac t\pi\right)^n\incn
e^{t\langle z,w\rangle}w_1^ke^{-(t+x)|w|^2}\,dv(w)\\
&=&\left(\frac t{t+x}\right)^n\incn e^{(t+x)\langle tz/(t
+x),w\rangle}w_1^k\,dv_{t+x}(w)\\
&=&\left(\frac t{t+x}\right)^n\left(\frac{tz_1}{t+x}
\right)^k\\
&=&\left(\frac t{t+x}\right)^{n+k}z_1^k.
\end{eqnarray*}

We next calculate the norm of $S_t(f_{x,k})$ in $L^p(\cn,dv_s)$
with the help of Lemma~\ref{1} again.
\begin{eqnarray*}
\incn|S_t(f_{x,k})|^p\,dv_s&=&\left(\frac t{t+x}\right)^{p
(n+k)}\incn|z_1|^{pk}\,dv_s(z)\\
&=&\left(\frac t{t+x}\right)^{p(n+k)}\frac{\Gamma((pk/2)+1)}
{s^{pk/2}}.
\end{eqnarray*}

Now if the integral operator $S_t$ is bounded on $L^p(\cn,dv_s)$,
then there exists a positive constant $C$ (independent of $x$ and $k$)
such that
$$\left(\frac t{t+x}\right)^{p(n+k)}\frac{\Gamma((pk/2)+1)}
{s^{pk/2}}\le C\left(\frac s{px+s}\right)^n
\frac{\Gamma((pk/2)+1)}{(px+s)^{pk/2}},$$
or
$$\left(\frac t{t+x}\right)^{p(n+k)}\le C\left(\frac{s}{s
+px}\right)^{n+(pk/2)}.$$
Fix any $x>0$ and look at what happens in the above inequality when
$k\to\infty$. We deduce that
$$\left(\frac t{t+x}\right)^2\le\frac s{s+px}.$$
Cross multiply the two sides of the inequality above and simplify.
The result is
$$pt^2\le 2st+sx.$$
Let $x\to0$. Then $pt^2\le 2st$, or $pt\le 2s$.
This completes the proof of the lemma.
\end{proof}

\begin{lemma}
Suppose $1<p<\infty$ and $S_t$ is bounded on $L^p(\cn,dv_s)$.
Then $pt>s$.
\label{7}
\end{lemma}

\begin{proof}
If $p>1$ and $S_t$ is bounded on $L^p(\cn,dv_s)$, then $S_t^*$ is
bounded on $L^q(\cn,dv_s)$, where $1/p+1/q=1$. Applying the formula for
$S_t^*$ from (\ref{eq4}) to the constant function $f=1$ shows that the
function $e^{(s-t)|z|^2}$ is in $L^q(\cn,dv_s)$. From this we
deduce that
$$q(s-t)<s,$$
which is easily seen to be equivalent to $s<pt$.
\end{proof}

\begin{lemma}
If $S_t$ is bounded on $L^1(\cn,dv_s)$, then $t=2s$.
\label{8}
\end{lemma}

\begin{proof}
Fix any $a\in\cn$ and consider the function
$$f_a(z)=\frac{e^{t\langle z,a\rangle}}{|e^{t\langle z,a\rangle}|},
\qquad z\in\cn.$$
Obviously, $\|f_a\|_\infty=1$ for every $a\in\cn$. On the other hand, it
follows from (\ref{eq4}) and Lemma~\ref{3} that
\begin{eqnarray*}
S_t^*(f_a)(a)&=&\left(\frac ts\right)^ne^{(s-t)|a|^2}\incn
|e^{t\langle w,a\rangle}|\,dv_s(w)\\
&=&\left(\frac ts\right)^ne^{(s-t)|a|^2}e^{t^2|a|^2/(4s)}.
\end{eqnarray*}
Since $S^*_t$ is bounded on $L^\infty(\cn)$, there exists a positive constant
$C$ such that
$$\left(\frac ts\right)^ne^{(s-t)|a|^2}e^{t^2|a|^2/(4s)}
\le\|S_t^*(f_a)\|_\infty\le C\|f_a\|_\infty=C$$
for all $a\in\cn$. This clearly implies that
$$s-t+\frac{t^2}{4s}\le0,$$
which is equivalent to
$$(2s-t)^2\le0.$$
Therefore, we have $t=2s$.
\end{proof}

\begin{lemma}
Suppose $1<p\le2$ and $S_t$ is bounded on $L^p(\cn,dv_s)$.
Then $pt=2s$.
\label{9}
\end{lemma}

\begin{proof}
Once again, we consider functions of the form
$$f_{x,k}(z)=e^{-x|z|^2}z_1^k,\qquad z\in\cn,$$
where $x>0$ and $k$ is a positive integer. It follows from (\ref{eq4}) and
Lemma~\ref{2} that
\begin{eqnarray*}
S_t^*(f_{x,k})(z)&=&\left(\frac t\pi\right)^ne^{(s-t)|z|^2}
\incn e^{t\langle z,w\rangle}w_1^ke^{-(s+x)|w|^2}\,dv(w)\\
&=&\left(\frac t{s+x}\right)^ne^{(s-t)|z|^2}\incn
e^{(s+x)\langle tz/(s+x),w\rangle}w_1^k\,dv_{s+x}(w)\\
&=&\left(\frac t{s+x}\right)^ne^{(s-t)|z|^2}
\left(\frac{tz_1}{s+x}\right)^k\\
&=&\left(\frac t{s+x}\right)^{n+k}e^{(s-t)|z|^2}z_1^k.
\end{eqnarray*}

Suppose $1<p\le2$ and $1/p+1/q=1$. If the operator $S_t$ is bounded on
$L^p(\cn,dv_s)$, then the operator $S_t^*$ is bounded on
$L^q(\cn,dv_s)$. So there exists a positive constant $C$, independent of
$x$ and $k$, such that
$$\incn|S_t^*(f_{x,k})|^q\,dv_s\le C\incn|f_{x,k}|^q\,dv_s.$$
It follows from the proof of Lemma~\ref{1} that
$$\incn|f_{x,k}|^q\,dv_s=\left(\frac s{qx+s}\right)^n
\frac{\Gamma((qk/2)+1)}{(qx+s)^{qk/2}}.$$
On the other hand, it follows from Lemma~\ref{7} and its proof that
$s-q(s-t)>0$, so the integral
$$I=\incn|S_t^*(f_{x,k})|^q\,dv_s$$
can be evaluated with the help of Lemma~\ref{1} as follows.
\begin{eqnarray*}
I&=&\left(\frac t{s+x}\right)^{
q(n+k)}\left(\frac s\pi\right)^n\incn|z_1|^{qk}e^{-(s-q(s-t))
|z|^2}\,dv(z)\\
&=&\left(\frac t{s+x}\right)^{q(n+k)}\left(\frac s{s-q(s-
t)}\right)^n\incn|z_1|^{qk}\,dv_{s-q(s-t)}(z)\\
&=&\left(\frac t{s+x}\right)^{q(n+k)}\left(\frac s{s-q(s-
t)}\right)^n\frac{\Gamma((qk/2)+1)}{(s-q(s-t))^{qk/2}}.
\end{eqnarray*}
Therefore,
$$\left(\frac t{s+x}\right)^{q(n+k)}\left(\frac s{s-q(s-
t)}\right)^n\frac{\Gamma((qk/2)+1)}{(s-q(s-t))^{qk/2}}$$
is less than or equal to
$$C\left(\frac s{qx+s}\right)^n
\frac{\Gamma((qk/2)+1)}{(qx+s)^{qk/2}},$$
which easily reduces to
$$\left(\frac t{s+x}\right)^{q(n+k)}\le C\left(\frac{s-
q(s-t)}{s+qx}\right)^{n+(qk/2)}.$$
Once again, fix $x>0$ and let $k\to\infty$. We find out that
$$\left(\frac t{s+x}\right)^2\le\frac{s-q(s-t)}{s+qx}.$$
Using the relation $1/p+1/q=1$, we can change the right-hand side above to
$$\frac{pt-s}{(p-1)s+px}.$$
It follows that
$$t^2(p-1)s+t^2px\le(pt-s)(s^2+2sx+x^2),$$
which can be written as
$$(pt-s)x^2+[2s(pt-s)-t^2p]x
+s^2(pt-s)-t^2(p-1)s\ge0.$$
Let $q(x)$ denote the quadratic function on the left-hand side of the above
inequality. Since $pt-s>0$ by Lemma~\ref{7}, the function $q(x)$
attains its minimum value at
$$x_0=\frac{pt^2-2s(pt-s)}{2(pt-s)}.$$
Since $2\ge p$, the numerator above is greater than or equal to
$$pt^2-2pts+ps^2=p(t-s)^2.$$
It follows that $x_0\ge0$ and so $h(x)\ge h(x_0)\ge0$ for all real $x$ (not
just nonnegative $x$). From this we deduce that the discriminant of $h(x)$
cannot be positive. Therefore,
$$[2s(pt-s)-pt^2]^2-4(pt-s)[s^2(pt-s)
-t^2(p-1)s]\le0.$$
Elementary calculations reveal that the above inequality is equivalent to
$$(pt-2s)^2\le0.$$
Therefore, $pt=2s$.
\end{proof}

\begin{lemma}
Suppose $2<p<\infty$ and $S_t$ is bounded on $L^p(\cn,dv_s)$.
Then $pt=2s$.
\label{10}
\end{lemma}

\begin{proof}
If $S_t$ is bounded on $L^p(\cn,dv_s)$, then $S_t^*$ is bounded
on $L^q(\cn,dv_s)$, where $1<q<2$ and $1/p+1/q=1$. It follows from
(\ref{eq4}) that there exists a positive constant $C$, independent
of $f$, such that
$$\incn\left|e^{(s-t)|z|^2}\incn e^{t\langle z,w\rangle}\left[f(w)
e^{(t-s)|w|^2}\right]\dva(w)\right|^q\,dv_s(z)$$
is less than or equal to
$$C\incn|f(w)|^q\,dv_s(w),$$
where $f$ is any function in $L^q(\cn,dv_s)$. Let
$$f(z)=g(z)e^{(s-t)|z|^2},$$
where $g\in L^q(\cn,dv_{s-q(s-t)})$ (recall from Lemma~\ref{7}
that $s-q(s-t)>0$). We obtain another positive constant C
(independent of $g$) such that
$$\incn|S_t g|^q\,dv_{s-q(s-t)}\le
C\incn|g|^q\,dv_{s-q(s-t)}$$
for all $g\in L^q(\cn,dv_{s-q(s-t)})$. Therefore, the operator
$S_t$ is bounded on $L^q(\cn,dv_{s-q(s-t)})$. Since $1<q<2$,
it follows from Lemma~\ref{9} that
$$qt=2[s-q(s-t)].$$
It is easy to check that this is equivalent to $pt=2s$.
\end{proof}

We now complete the proof of the first part of the main theorem. As was 
pointed out in the introduction, this part of the theorem is known before.
We included a full proof here for two purposes. First, this gives a 
different and self-contained approach. Second, as a by-product of this
different approach, we are going to obtain the inequality $\|T_t\|\le 2^n$,
which is one half of the identity $\|T_t\|=2^n$.

\begin{thm}
Suppose $t>0$, $s>0$, and $p\ge1$. Then the following conditions
are equivalent.
\begin{enumerate}
\item[(a)] The operator $T_t$ is bounded on $L^p(\cn,dv_s)$.
\item[(b)] The operator $S_t$ is bounded on $L^p(\cn,dv_s)$.
\item[(c)] The weight parameters satisfy $pt=2s$.
\end{enumerate}
\label{11}
\end{thm}

\begin{proof}
When $p=1$, that (b) implies (c) follows from Lemma~\ref{8},
that (c) implies (a) follows from Fubini's theorem and Lemma~\ref{3}, and
that (a) implies (b) is obvious.

When $1<p<\infty$, that (b) implies (c) follows from Lemmas~\ref{9}
and \ref{10}, and that (a) implies (b) is still obvious.

So we assume $1<p<\infty$ and proceed to show that condition (c) implies (a).
We do this with the help of Schur's test (Lemma~\ref{5}).

Let $1/p+1/q=1$ and consider the positive function
$$h(z)=e^{\lambda|z|^2},\qquad z\in\cn,$$
where $\lambda$ is a constant to be specified later.

Recall that
$$T_t f(z)=\incn H(z,w)f(w)\,dv_s(w),$$
where
$$H(z,w)=\left(\frac ts\right)^n|e^{t\langle z,w\rangle}
e^{(s-t)|w|^2}|$$
is a positive kernel. We first consider the integrals
$$I(z)=\incn H(z,w)h(w)^q\,dv_s(w),\qquad z\in\cn.$$
If $\lambda$ satisfies
\begin{equation}
t>q\lambda,
\label{eq6}
\end{equation}
then it follows from Lemma~\ref{3} that
\begin{eqnarray*}
I(z)&=&\left(\frac t\pi\right)^n\incn|e^{t\langle z,w\rangle}|\,e^{-(t-
q\lambda)|w|^2}\,dv(w)\\
&=&\left(\frac t{t-q\lambda}\right)^n\incn|e^{t\langle z,w\rangle}|
\,dv_{t-q\lambda}(w)\\
&=&\left(\frac t{t-q\lambda}\right)^ne^{t^2|z|^2/4(t-q\lambda)}.
\end{eqnarray*}
If we choose $\lambda$ so that
\begin{equation}
\frac{t^2}{4(t-q\lambda)}=q\lambda,
\label{eq7}
\end{equation}
then we obtain
\begin{equation}
\incn H(z,w)h(w)^q\,dv_s(w)\le\left(\frac t{t-q\lambda}\right)^nh(z)^q
\label{eq8}
\end{equation}
for all $z\in\cn$.

We now consider the integrals
$$J(w)=\incn H(z,w)h(z)^p\,dv_s(z),\qquad w\in\cn.$$
If $\lambda$ satisfies
\begin{equation}
s-p\lambda>0,
\label{eq9}
\end{equation}
then it follows from Lemma~\ref{3} that
\begin{eqnarray*}
J(w)&=&\left(\frac ts\right)^n\incn|e^{t\langle z,w\rangle}e^{(s-t)
|w|^2}|\,h(z)^p\,dv_s(z)\\
&=&\left(\frac t\pi\right)^ne^{(s-t)|w|^2}\incn|e^{t\langle z,w
\rangle}|\,e^{-(s-p\lambda)|z|^2}\,dv(z)\\
&=&\left(\frac t{s-p\lambda}\right)^ne^{(s-t)|w|^2}
e^{t^2|w|^2/4(s-p\lambda)}\\
&=&\left(\frac t{s-p\lambda}\right)^n
e^{[(s-t)+t^2/4(s-p\lambda)]|w|^2}.
\end{eqnarray*}
If we choose $\lambda$ so that
\begin{equation}
s-t+\frac{t^2}{4(s-p\lambda)}=p\lambda,
\label{eq10}
\end{equation}
then we obtain
\begin{equation}
\incn H(z,w)h(z)^p\,dv_s(z)\le\left(\frac t{s-p\lambda}\right)^nh(w)^p
\label{eq11}
\end{equation}
for all $w\in\cn$. In view of Schur's test and the estimates in (\ref{eq8})
and (\ref{eq11}), we conclude that the operator $T_t$ would be bounded
on $L^p(\cn,dv_s)$ provided that we could choose a real $\lambda$ to satisfy
conditions (\ref{eq6}), (\ref{eq7}), (\ref{eq9}), and (\ref{eq10})
simultaneously.

Under our assumption that $pt=2s$ it is easy to verify that condition (\ref{eq7})
is the same as condition (\ref{eq10}). In fact, we can explicitly solve for $q\lambda$
and $p\lambda$ in (\ref{eq7}) and (\ref{eq10}), repectively, to obtain
$$q\lambda=\frac t2,\quad p\lambda=\frac{2s-t}2.$$
The relations $pt=2s$ and $1/p+1/q=1$ clearly imply that the two resulting $\lambda$'s
above are consistent, namely,
\begin{equation}
\lambda=\frac t{2q}=\frac{2s-t}{2p}.
\label{eq12}
\end{equation}
Also, it is easy to see that the above choice of $\lambda$ satisfies both (\ref{eq6})
and (\ref{eq9}). This completes the proof of the theorem.
\end{proof}

\begin{thm}
If $1\le p<\infty$ and $pt=2s$, then
$$\incn|S_t f|^p\,dv_s\le\incn|T_t f|^p\,dv_s\le
2^{np}\incn|f|^p\,dv_s$$
for all $f\in L^p(\cn,dv_s)$.
\label{12}
\end{thm}

\begin{proof}
With the choice of $\lambda$ in (\ref{eq12}), the constants in (\ref{eq8})
and (\ref{eq11}) both reduce to $2^n$. Therefore, Schur's test tells us that,
in the case when $1<p<\infty$, the norm of $T_t$ on $L^p(\cn,dv_s)$ does not 
exceed $2^n$.

When $p=1$, the desired estimate follows from Fubini's theorem and Lemma~\ref{3}.
\end{proof}

Theorem~\ref{12} above can be stated as $\|S_t\|\le\|T_t\|\le2^n$, with $S_t$
and $T_t$ considered as operators on $L^p(\cn,dv_s)$. We now proceed to the proof 
of the inequality $\|T_t\|\ge2^n$. Several lemmas are needed for this estimate.

\begin{lemma}
For $c>0$ and $p\geq 1$ we have
$$\lim_{h\to0^{+}}h\!\int_{c}^{\infty }\!\left[\int_{c}^{\infty }(uv)^{-\frac{1}{4}}
\exp\left(\!\sqrt{uv}-\frac{u+v}{2}-\frac{hv}{p}\right)dv\right]^pdu=
\left(2\sqrt{2\pi}\right)^{p}.$$
\label{13}
\end{lemma}

\begin{proof}
We begin with the inner integral
$$I(u)=\int_c^\infty(uv)^{-\frac14}\exp\left(\sqrt{uv}-\frac{u+v}2-\frac{hv}p
\right)\,dv.$$
Let $a=\frac{1}{2}+\frac{h}{p}$ and change variables according to $v=t^2$. Then
$$I(u)=2u^{-\frac14}e^{-\frac u2}\int_{\sqrt c}^\infty\sqrt t\exp(-at^2+
\sqrt u\,t)\,dt.$$
Write
$$-at^2+\sqrt u\,t=-a\left(t-\frac{\sqrt u}{2a}\right)^2+\frac u{4a},$$
make a change of variables according to $x=t-(\sqrt u/2a)$, and simplify the
result. We obtain
$$I(u)=2u^{-\frac14}e^{-\frac{uh}{2ap}}\int_{\sqrt c-\frac{\sqrt u}{2a}}^\infty
\sqrt{x+\frac{\sqrt u}{2a}}\,e^{-ax^2}\,dx.$$
It is then clear that we can rewrite $I(u)$ as follows.
$$I(u)=\varphi_1(u)+\varphi_2(u),$$
where
$$\varphi_1(u)=\frac{2}{\sqrt{2a}}\,e^{-\frac{uh}{2ap}}
\int_{\sqrt{c}-\frac{\sqrt{u}}{2a}}^{+\infty}e^{-ax^{2}}\,dx,$$
and
$$\varphi_2(u)=\frac{2\,e^{-\frac{uh}{2ap}}}{u^{\frac{1}{4}}}
\int_{\sqrt{c}-\frac{\sqrt{u}}{2a}}^{+\infty}e^{-ax^{2}}\left(\sqrt{x+
\frac{\sqrt{u}}{2a}}-\sqrt{\frac{\sqrt{u}}{2a}}\,\right)\,dx.$$

For the function $\varphi_2$ we rationalize the numerator in its integrand to obtain
\begin{eqnarray*}
\left\vert\varphi_2(u)\right\vert &\leq &\frac{2}{u^{\frac{1}{4}}}\,
e^{-\frac{uh}{2ap}}\int_{\sqrt{c}-\frac{\sqrt{u}}{2a}}^{+\infty
}e^{-ax^{2}}\frac{|x|}{\sqrt{\frac{\sqrt{u}}{2a}}}\,dx\\
&\leq &\frac{2\sqrt{2a}}{u^{\frac{1}{2}}}\,e^{-\frac{uh}{2ap}}
\int_{-\infty}^{+\infty }|x|\,e^{-ax^{2}}\,dx.
\end{eqnarray*}
A simple calculation of the last integral above then gives
\begin{equation}
\left\vert\varphi_2(u)\right\vert \leq \frac{4}{\sqrt{2a}}
\,u^{-\frac{1}{2}}\,e^{-\frac{uh}{2ap}}.  \label{meq14}
\end{equation}
Similarly, we have
\begin{equation}
\left\vert \varphi_1(u)\right\vert \leq \frac2{\sqrt{2a}}e^{-\frac{uh}{2ap}}
\int_{-\infty}^\infty e^{-ax^2}\,dx=\frac{\sqrt{2\pi}}{a}\,e^{-\frac{uh}{2ap}}.
\label{meq15}
\end{equation}

We now use the above estimates to show that
\begin{equation}
\lim_{h\to 0^{+}}h\int_{c}^{\infty}\varphi_1(u)^{p-1}
\left\vert\varphi_2(u)\right\vert\,du=0, \label{meq16}
\end{equation}
and
\begin{equation}
\lim_{h\to 0^{+}}h\int_{c}^{\infty}\left\vert\varphi_2(u)
\right\vert^{p}\,du=0.
\label{meq17}
\end{equation}%
In fact, according to (\ref{meq14}) and (\ref{meq15}),
$$h\varphi_1(u)^{p-1}\left\vert\varphi_2(u)\right\vert 
\leq \left(\frac{\sqrt{2\pi}}{a}\right)^{p-1}\frac{4}
{\sqrt{2a}}\,h\,u^{-\frac{1}{2}}\,e^{-\frac{uh}{2a}},$$
from which we derive that
\begin{eqnarray*}
h\int_{c}^{\infty }\varphi_1(u)^{p-1}\left\vert
\varphi_2(u)\right\vert\,du &\leq &\left(\frac{\sqrt{2\pi }}
{a}\right)^{p-1}\,h\,\frac{4}{\sqrt{2a}}\int_{c}^{\infty }
u^{-\frac{1}{2}}\,e^{-\frac{uh}{2a}}\,du\\
&=&4\left(\frac{\sqrt{2\pi }}{a}\right)^{p-1}\sqrt{h}
\int_{\frac{ch}{2a}}^{+\infty }\omega^{-\frac{1}{2}}e^{-\omega}\,d\omega.
\end{eqnarray*}
Since $a\to\frac{1}{2}$ as $h\to0^{+}$, we obtain (\ref{meq16}).
On the other hand, it follows from (\ref{meq14}) that
$$\left\vert\varphi_2(u)\right\vert^{p}\leq \left(\frac{4}{\sqrt{2a}}
\right)^{p}u^{-\frac{p}{2}}\,e^{-\frac{uh}{2a}},$$
so
\begin{eqnarray*}
h\int_{c}^{\infty}\left\vert\varphi_2(u)\right\vert^{p}\,du &\leq &
\left(\frac{4}{\sqrt{2a}}\right)^{p}h\int_{c}^{\infty}u^{-\frac{p}{2}}\,
e^{-\frac{uh}{2a}}\,du=\\
&=& 2^{p-1}a^{1-p}h^{p/2}\int_{\frac{ch}{2a}}^{+\infty }
\omega^{-\frac{p}{2}}e^{-\omega}\,d\omega.
\end{eqnarray*}%
Let $h\to0^+$ and use the fact that $a\to\frac{1}{2}$ as $h\to0^{+}$.
We obtain (\ref{meq17}).

By the change of variables $s=uh/2a$, we have
\begin{eqnarray*}
h\int_{c}^{\infty}\varphi_1(u)^{p}\,du&=&h\int_{c}^{\infty}\left(\frac{2}
{\sqrt{2a}}\right)^{p}e^{-\frac{uh}{2a}}\left(\int_{\sqrt{c}-
\frac{\sqrt{u}}{2a}}^{+\infty }e^{-ax^{2}}dx\right)^{p}du\\
&=&\frac{2^{p+1}a}{\left(\sqrt{2a}\right)^{p}}\int_{\frac{hc}{2a}}^{\infty}
e^{-s}\left(\int_{\sqrt{c}-\frac{\sqrt{s}}{\sqrt{2ah}}}^{+\infty}e^{-ax^{2}}
dx\right)^{p}ds.
\end{eqnarray*}
Let $h\to0^+$, notice that $a\to\frac12$ as $h\to0^+$, and use Lebesgue's 
dominated convergence theorem. We get
\begin{equation}
\lim_{h\to0^{+}}h\!\int_{c}^{\infty}\!\varphi_1^{p}\,du=2^{p}\!\int_{0}^{\infty}
e^{-s}\left[\int_{-\infty }^{+\infty}e^{-\frac{x^{2}}{2}}dx\right]^{p}ds=
\left(2\sqrt{2\pi}\right)^{p}.  \label{meq18}
\end{equation}

If $p\ge1$, it is easy to see that the function
$$g(z)=\frac{|1+z|^p-|z|^p}{1+|z|^{p-1}}$$
is continuous and bounded on $\C$. Replacing $z$ by $z/w$, we see that
$$\bigl||z+w|^p-|z|^p\bigr|\le C\bigl(|z|^{p-1}|w|+|w|^p\bigr)$$
for all $z$ and $w$, where $C$ is a positive constant that only
depends on $p$. This along with (\ref{meq16}) and (\ref{meq17}) shows that
$$\lim_{h\to0^{+}}h\int_{c}^{\infty}\left\vert\left(\varphi_1(u)+\varphi_2(u)
\right)^{p}-\varphi_1(u)^{p}\right\vert\,du=0.$$
Combining this with (\ref{meq18}), we conclude that
$$\lim_{h\to0^+}h\int_c^\infty\left(\varphi_1(u)+\varphi_2(u)\right)^p\,du=
\left(2\sqrt{2\pi}\,\right)^p.$$
This proves the desired result.
\end{proof}

\begin{lemma}
Let 
$$\mathcal{K}(x,y)=\sum_{n=0}^{\infty}\frac{t^{2n}x^{n}y^{n}}{4^{n}(n!)^{2}}$$
and define an integral operator
$$A:L^p(0,\infty)\to L^p(0,\infty)$$
by
$$Af(x)=\int_{0}^{\infty }te^{-\frac{t}{2}\left(x+y\right)}\mathcal{K}(x,y)
f(y)\,dy,$$
where $p\ge1$ and $t$ is any fixed positive constant. Then the norm of $A$ on
$L^p(0,\infty)$ satisfies $\left\Vert A\right\Vert \geq 2$.
\label{14}
\end{lemma}

\begin{proof}
It follows from the asymptotic behavior of the Bessel function $J_0(x)$ (see
page 199 of \cite{W} for example) that
$$\sum_{k=0}^{\infty}\frac{(z/2)^{2k}}{\left(k!\right)^{2}}\sim
\frac{e^{z}}{\sqrt{2\pi z}}$$
as $z\to\infty$. Thus%
$$\sum_{n=0}^{\infty}\frac{u^{n}}{\left(n!\right)^{2}}
\sim\frac{e^{2\sqrt{u}}}{\sqrt{4\pi}\,u^{\frac{1}{4}}}$$
as $u\to\infty$. Fix an arbitrary $\eta>0$ and choose some $u_{0}>0$ such that
\begin{equation}
\sum_{n=0}^{\infty}\frac{u^{n}}{\left(n!\right)^{2}}>\left(1-\eta\right)
\frac{e^{^{2\sqrt{u}}}}{\sqrt{4\pi}\,u^{\frac{1}{4}}}  \label{meq19}
\end{equation}%
for every $u\geq u_{0}$.

Let $\delta_{0}=2\sqrt{u_{0}}/t$. It follows from (\ref{meq19}) that
\begin{equation}
\mathcal{K}(x,y)>\left(1-\eta\right)\frac{e^{^{t\sqrt{xy}}}}
{\left(xy\right)^{\frac{1}{4}}\sqrt{2\pi t}}  \label{meq20}
\end{equation}
for all $x\ge\delta_0$ and all $y\ge\delta_0$.
Fix a positive number $\varepsilon$ and let 
$$f_{\varepsilon}(x)=\exp\left(-\frac{\varepsilon x}{p}\right).$$
Then 
$$\left\Vert f_{\varepsilon}\right\Vert=\left[\int_{0}^{\infty}
\left\vert f_{\varepsilon}(x)\right\vert^{p}\,dx\right]^{1/p}=
\varepsilon^{-\frac{1}{p}},$$
and so
\begin{equation}
\left\Vert A\right\Vert \geq \varepsilon ^{\frac{1}{p}}\left\Vert
Af_{\varepsilon}\right\Vert.  \label{meq21}
\end{equation}
On the other hand, it follows from (\ref{meq20}) that
\begin{eqnarray*}
\left\Vert Af_{\varepsilon}\right\Vert &\geq &\left(\int_{\delta_{0}}^{\infty}
\left\vert Af_{\varepsilon}(x)\right\vert^{p}\,dx\right)^{\frac{1}{p}}\\
& \geq & \left[\int_{\delta_{0}}^{\infty}\,dx\left(\,\int_{\delta_{0}}^{\infty}
te^{-\frac{t}{2}(x+y)}\mathcal{K}(x,y)e^{-\frac{\varepsilon y}{p}}\,dy
\right)^{p}\,\right]^{\frac{1}{p}} \\
&\geq &\frac{(1-\eta)t}{\sqrt{2\pi t}}\left[\int_{\delta _{0}}^{\infty }
\,dx\left(\,\int_{\delta_{0}}^{\infty }\frac{e^{-\frac{t}{2}\left(x+y\right)+
t\sqrt{xy}-\frac{\varepsilon y}{p}}}{\left( xy\right)^{\frac{1}{4}}}\,dy\right)^{p}
\,\right]^{\frac{1}{p}}.
\end{eqnarray*}
Combining this with (\ref{meq21}), we obtain
$$\left\Vert A\right\Vert \geq \frac{(1-\eta)t}{\sqrt{2\pi t}}
\left[\varepsilon\int_{\delta_{0}}^{\infty}\,dx\left(\,
\int_{\delta _{0}}^{\infty}\frac{e^{-\frac{t}{2}(x+y)+t
\sqrt{xy}-\frac{\varepsilon y}{p}}}{\left( xy\right) ^{\frac{1}{4}}}
\,dy\right)^{p}\,\right]^{\frac{1}{p}}.$$
After the change of variables $xt=u$ and $yt=v$ we obtain
\begin{equation}
\left\Vert A\right\Vert\geq\frac{1-\eta}{\sqrt{2\pi}}
\left[\frac{\varepsilon}{t}\int_{t\delta_{0}}^{\infty }du\left(\,\,
\int_{t\delta_{0}}^{\infty}\frac{e^{-\frac{1}{2}(u+v)+\sqrt{uv}-\frac{\varepsilon v}
{tp}}}{(uv)^{\frac{1}{4}}}\,dv\right)^{p}\,\right]^{\frac{1}{p}}.  \label{meq22}
\end{equation}%

Let $t\delta_{0}=c$ and $\varepsilon/t=h$. Clearly, $h\to 0^{+}$ when $\varepsilon
\to 0^{+}$. Let $\varepsilon\to0^+$ in (\ref{meq22}) and apply Lemma~\ref{13}. We obtain
$$\left\Vert A\right\Vert \geq \frac{1-\eta}{\sqrt{2\pi}}\cdot 2\sqrt{2\pi}=2(1-\eta).$$
Since $\eta >0$ is arbitrary, we obtain $\left\Vert A\right\Vert \geq 2$, and the
proof of Lemma~\ref{14} is complete.
\end{proof}

We can now prove the main result of the paper.

\begin{thm}
If $1\leq p<\infty$ and $pt=2s$, then the norm of $\,T_t$ on $L^p(\cn,dv_s)$ is given
by $\left\Vert T_{t}\right\Vert=2^{n}$.
\label{15}
\end{thm}

\begin{proof}
In view of Theorem~\ref{12} it is enough for us to prove the inequality $\|T_t\|\ge2^n$.

Recall that when $n=1$, we use the notation $dA_s$ instead of $dv_s$. For 
$f\in L^p(\C,dA_s)$ we consider
$$\Phi(z_1,\cdots,z_n)=f(z_{1})\cdots f(z_{n}).$$
Then $\Phi\in L^p(\cn,dv_s)$ and we have
$$\left\Vert T_{t}\right\Vert^{p}\geq\frac{\left\Vert T_{t}\Phi\right\Vert^{p}}
{\left\Vert\Phi\right\Vert^{p}}=\left[\frac{\displaystyle\int_{\C}\left|
\int_{\C}|e^{tz\,\overline\zeta}|f(\zeta)\,dA_{t}(\zeta)
\right|^{p}dA_s(z)}{\displaystyle\int_{\C}\left\vert f(\zeta)\right\vert^{p}
\,dA_{s}(\zeta)}\right]^n.$$
When $f$ runs over all unit vectors in $L^p(\C,dA_s)$, the supremum of the quotient 
inside the brackets above is exactly the $p$th power of the norm of the operator 
$T_t$ on $L^p(\C,dA_s)$. So we only need to prove the inequality $\|T_t\|\ge2^n$ 
for $n=1$.

Now we assume $n=1$, $p\ge1$, and let
$$T_t:L^p(\C,dA_s)\to L^p(\C,dA_s)$$
be the integral operator defined by
$$T_{t}f(z)=\int_{\C}|e^{tz\,\overline\zeta}|\,f(\zeta)\,dA_t(\zeta).$$
To obtain a lower estimate of the norm of $T_t$ on $L^p(\C,dA_s)$, we apply $T_t$ to
a family of special functions. More specifically, we consider functions of the form
$$f(z)=G(|z|^{2})e^{t|z|^2/2},\qquad z\in\C,$$
where $G$ is any unit vector in $L^p(0,\infty)$. It follows from polar coordinates 
and the assumption $pt=2s$ that
\begin{equation}
\|f\|^p=\int_\C|f|^p\,dA_s=s\int_0^\infty|G(x)|^p\,dx=s.\label{meq24}
\end{equation}%
On the other hand,
\begin{eqnarray*}
T_tf(z)&=&\int_\C|e^{tz\overline w}|G(|w|^2)e^{t|w|^2/2}\,dA_t(w)\\
&=&\int_\C|e^{tz\overline w/2}|^2G(|w|^2)e^{t|w|^2/2}\,dA_t(w)\\
&=&\int_\C\left|\sum_{n=0}^\infty\frac{(tz\overline w/2)^n}{n!}\right|^2
G(|w|^2)e^{t|w|^2/2}\,dA_t(w)\\
&=&\sum_{n=0}^\infty\frac{t^{2n}|z|^{2n}}{4^n(n!)^2}\int_\C|w|^{2n}
G(|w|^2)e^{t|w|^2/2}\,dA_t(w)\\
&=&\sum_{n=0}^\infty\frac{t^{2n}|z|^{2n}}{4^n(n!)^2}\,t\int_0^\infty
y^nG(y)e^{-ty/2}\,dy\\
&=&\int_0^\infty te^{-ty/2}\mathcal{K}(|z|^2,y)G(y)\,dy\\
&=&e^{t|z|^2/2}AG(|z|^2),
\end{eqnarray*}
where the kernel $\mathcal{K}$ and the operator $A$ are from Lemma~\ref{14}.
Using polar coordinates and the assumption $pt=2s$ one more time, we obtain
\begin{equation}
\|T_tf\|^p=\int_\C|T_tf(z)|^p\,dA_s(z)=s\int_0^\infty|AG(x)|^p\,dx.
\label{meq25}
\end{equation}
By (\ref{meq24}) and (\ref{meq25}) we must have
\begin{equation}
\left\Vert T_{t}\right\Vert^{p}\geq \frac{\left\Vert T_{t}f\right\Vert^{p}}
{\left\Vert f\right\Vert^{p}}=\int_0^\infty|AG(x)|^p\,dx.  \label{meq26}
\end{equation}%
Take the supremum over $G$ and apply Lemma~\ref{14}. The result is
$$\left\Vert T_{t}\right\Vert^{p}\geq \left\Vert A\right\Vert^{p}\geq2^{p}.$$
This completes the proof of the theorem.
\end{proof}

We conclude the paper with two corollaries.

\begin{cor}
For any $s>0$ and $p\ge1$ the Fock space $F^p_s$ is a complemented subspace of
$L^p(\cn,dv_s)$, that is, there exists a closed subspace $X^p_s$ of $L^p(\cn,dv_s)$
such that
$$L^p(\cn,dv_s)=F^p_s\oplus X^p_s,$$
where $\oplus$ denotes the direct sum of two subspaces.
\label{16}
\end{cor}

\begin{proof}
Choose $t>0$ such that $pt=2s$. Then by Theorem~\ref{11}, the operator $S_t$ is
a bounded projection from $L^p(\cn,dv_s)$ onto $F^p_s$. This shows that $F^p_s$
is complemented in $L^p(\cn,dv_s)$.
\end{proof}

The following result is obviously a generalization of Theorem~\ref{11}, but it is
also a direct consequence of Theorem~\ref{11}.

\begin{cor}
Suppose $a>0$, $b>0$, $s>0$, and $p\ge1$. Then the following conditions are
equivalent.
\begin{enumerate}
\item[(a)] The integral operator
$$T_{a,b}f(z)=\incn|e^{-a|z|^2+(a+b)\langle z,w\rangle-b|w|^2}|\,f(w)\,dv(w)$$
is bounded on $L^p(\cn,dv_s)$.
\item[(b)] The integral operator
$$S_{a,b}f(z)=\incn e^{-a|z|^2+(a+b)\langle z,w\rangle-b|w|^2}\,f(w)\,dv(w)$$
is bounded on $L^p(\cn,dv_s)$.
\item[(c)] The parameters satisfy $p(a+b)=2(s+pa)$.
\end{enumerate}
\label{17}
\end{cor}

\begin{proof}
The boundedness of $S_{a,b}$ on $L^p(\cn,dv_s)$ is equivalent to the existence of
a positive constant $C$, independent of $f$, such that
$$\incn\left|\incn e^{(a+b)\langle z,w\rangle-b|w|^2}\,f(w)\,dv(w)\right|^p\,dv_{s+pa}(z)$$
is less than or equal to
$$C\incn|f(z)e^{a|z|^2}|^p\,dv_{s+pa}(z).$$
Replacing $f(z)$ by $f(z)e^{-a|z|^2}$, we see that the above condition is equivalent to
$$\incn\left|\incn e^{(a+b)\langle z,w\rangle}\,f(w)\,dv_{a+b}(w)\right|^p\,dv_{s+pa}(z)
\le C\incn|f|^p\,dv_{s+pa}.$$
This is clearly equivalent to the boundedness of $S_{a+b}$ on $L^p(\cn,dv_{s+pa})$,
which, according to Theorem~\ref{11}, is equivalent to $p(a+b)=2(s+pa)$. Therefore,
conditions (b) and (c) are equivalent. The equivalence of (a) and (c) is proved
in exactly the same way.
\end{proof}

\end{document}